# THE TENSION BETWEEN INTUITIVE INFINITESIMALS AND FORMAL MATHEMATICAL ANALYSIS

## Mikhail G. Katz & David Tall

Abstract: we discuss the repercussions of the development of infinitesimal calculus into modern analysis, beginning with viewpoints expressed in the nineteenth and twentieth centuries and relating them to the natural cognitive development of mathematical thinking and imaginative visual interpretations of axiomatic proof.

## 1. Klein's reflections on "mystical schemes" in the calculus

Infinitesimal calculus is a dead metaphor. In countless courses of instruction around the globe, students register for courses in "infinitesimal calculus" only to find themselves being trained to perform epsilontic multiple-quantifier logical stunts, or else being told briefly about "the rigorous approach" to limits, promptly followed by instructions not to worry about it.

Anticipating the problem as early as 1908, Felix Klein reflected upon the success of a calculus textbook dealing in "mystical schemes", namely

> the textbook by Lübsen [...] which appeared first in 1855 and which had for a long time an extraordinary influence among a large part of the public [...] Lübsen defined the differential quotient first by means of the limit notion; but along side of this he placed [...] what he considered to be the *true infinitesimal calculus*—a mystical scheme of operating with infinitely small quantities [...] And then follows an English quotation: "An infinitesimal is the spirit of a departed quantity" [17, p. 216-217].

In his visionary way, Klein adds:

> The reason why such reflections could so long hold their place [alongside] the mathematically rigorous method of limits, must be sought probably in the widely felt need of penetrating beyond the abstract logical formulation of the method of limits to the intrinsic nature of continuous magnitudes, and of forming more definite images of them than were supplied by emphasis solely upon the psychological moment which determined the concept of limit [17, p. 217].

## 2. Interesting infinitesimals lead to contradictions

In the closing months of World War II, the teenage Peter Roquette's calculus teacher at Königsberg was an old lady trained in the old school, the regular teacher having been drafted into action. Roquette reminisces in the following terms:

> I still remember the sight of her standing in front of the blackboard w[h]ere she had drawn a wonderfully smooth parabola, inserting a secant and telling us that $\Delta y/\Delta x$ is its slope, until finally she convinced us that the slope of the tangent is *dy/dx* where *dx* is infinitesimally small and *dy* accordingly [26, p. 186].

Roquette recalls his youthful reaction:

> This, I admit, impressed me deeply. Until then our school Math had consisted largely of Euclidean geometry, with so many problems of constructing triangles from some given data. This was o.k. but in the long run that stuff did not strike me as more than boring exercises. But now, with those infinitesimals, Math seemed to have more interesting things in stock than I had met so far [26, p. 186].

But then at the university a few years later,

> we were told to my disappointment that my Math teacher had not been up to date after all. We were warned to beware of infinitesimals since they do not exist, and in any case they lead to contradictions. Instead, although one writes *dy/dx* […], this does not really mean a quotient of two entities, but it should be interpreted as a symbolic notation only, namely the limit of the quotient $\Delta y/\Delta x$. I survived this disappointment too [26, p. 186-187].

Then, some decades later, the old lady turned out not to have been so far off the mark:

> when I learned about Robinson's infinitesimals [24], my early school day experiences came to my mind again and I wondered whether that lady teacher had not been so wrong after all. The discussion with Abraham Robinson kindled my interest and I wished to know more about it. Some time later there arose the opportunity to invite him to visit us in Germany where he gave lectures on his ideas, first in Tübingen and later in Heidelberg, after I had moved there [26, p. 187].

The results of the ensuing collaboration were reported in [25] and [27].

Roquette mentions an infinitesimal calculus textbook published as late as 1912, the year of the last edition of L. Kiepert [16]. He speculates [26, p. 192] that his old lady teacher may have been trained using Kiepert's textbook.

## 3. Courant and infinitesimals "devoid of meaning"

Kiepert and other infinitesimal textbooks seem to have been edged out of the market by Courant's textbook [6]. Courant set the tone for the attitude prevailing at the time, when he described infinitesimals as "devoid of any clear meaning" and "naive befogging" [6, p. 81], as well as "incompatible with the clarity of ideas demanded in mathematics", "entirely meaningless", "fog which hung round the foundations", and a "hazy idea" [6, p. 101], while acknowledging Leibniz's masterly use of them:

> In the early days of the differential calculus even Leibnitz[1] himself was capable of combining these vague mystical ideas with a thoroughly clear understanding of the limiting process. It is true that this fog which hung round the foundations of the new science did not prevent Leibnitz or his great successors from finding the right path [6, p. 101].

How is it that they were in a position to find the right path? The Russian mathematician and historian Medvedev asks the million dollar question:

> If infinitely small and infinitely large magnitudes are regarded as inconsistent notions, how could they serve as a basis for the construction of so [magnificent] an edifice of one of the most important mathematical disciplines? [20], [21].

## 4. Vygodskiĭ: from biped back to quadruped?

In a 1931 letter [19] to the mathematician Vygodskiĭ, Luzin presents a hilarious account of the reception of Vygodskiĭ's infinitesimal calculus textbook in Soviet Russia. Vygodskiĭ dared to exploit *actual* infinitesimals. Luzin describes the reactions that ensued, in the following terms:

> I heard talk in Moscow about the restoration of the phlogiston[2] theory in science and charges of decadence [19, p. 68].

A modern reader may need to be reminded that in Stalinist Russia, a charge of bourgeois decadence was not to be trifled with, and could lead to a lengthy term in Siberian bestiaria or worse. The defenders of ideological (and decidedly secular) purity did not stop at invocations of phlogiston:

> In Leningrad [...] I heard talk to the effect that while Darwin [traced] the path of man's evolution from quadruped to biped, efforts are underway in mathematics to reverse this course [19, p. 68].

In a show of solidarity with Vygodskiĭ, Luzin proceeds to endorse a viewpoint strikingly similar to Cauchy's 1821 text [4] (which was apparently unavailable to Luzin):

> Unlike my colleagues, I think that an attempt to reconsider the idea of an infinitesimal as a variable finite quantity is fully scientific, and that the proposal to replace variable infinitesimals by fixed ones, far from having purely pedagogical significance, has in its favor something immeasurably deeper, and that this idea is growing roots in modern analysis [19, p. 68].

Luzin notes that

> the idea of the actually infinitely small has certain deep roots in the mind [19, p. 68].[3]

---

[1] In English speaking countries, the German name 'Leibniz' is often transliterated to 'Leibnitz' to represent the sound rather than the original spelling.
[2] Phlogiston was once thought to be a fire-like element contained within combustible bodies, and released during combustion, but became incommensurable 250 years ago.
[3] Lakoff & Núñez [18] would certainly agree.

In a possible allusion to despotic pre-revolutionary Russia of his student years, Luzin notes:

> The theory of limits entered my mind mechanically and crudely, not in a refined way but rather in a forced, police-like manner [19, p. 70].

The stark choice between Weierstrassian limits and infinitesimals came in Luzin's sophomore year:

> When the professors announced that *dy/dx* is the limit of a ratio, I thought: "What a bore! Strange and incomprehensible. No! They won't fool me: it's simply the ratio of infinitesimals, nothing else." [19, p. 70.]

Luzin's appropriately sophomoric attempt to construct a simpler version of Weierstrass's nowhere differentiable curve, by means of a diagonal "saw" (with numerous "steps" climbing along the diagonal of a square, see Figure 1) with infinitely many infinitesimal teeth, was patiently rebuffed by Professor Boleslav Kornelievich Mlodzeevskiĭ[4] (whom we later refer to as M.) on the grounds that "the actually infinite does not exist".

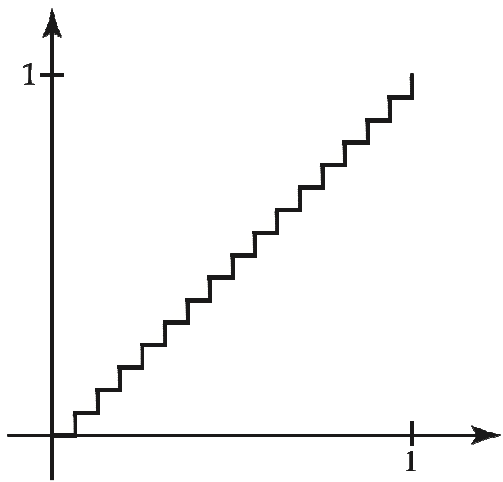

*Figure 1. Luzin's saw*

Unfazed, Luzin pursued M., after M.'s lecture on Cantor's cardinalities and $\aleph_0$. Luzin thought that

> These are complete contradictions: in analysis they say that every number is finite and modestly pass over in silence points at infinity on straight lines. In geometry, on the contrary, they keep on talking about points at infinity and deduce marvelous things.[5] A week ago, Boleslav Kornelievich cut me short by explaining that the actually infinite does not exist'. And now he does it himself! [19, p. 72.]

---

[4] Mlodzeevskiĭ, who brought the ideas of Hilbert and Klein from Göttingen to Moscow and was the first professor in Moscow to lecture on set theory and the theory of functions.
[5] Possibly an allusion to projective geometry.

## 5. Luzin and his infinitesimal saw

Encouraged by his Cantorian insight, Luzin confronted M., this time armed with a diagonal saw with what he claimed were countably many teeth. M., patiently, countered with the claim that Luzin's saw is merely

> verbal but not real... it is not genuine. [19, p. 73.]

Luzin countered by asking whether

> the Weierstrass curve, is it genuine or logical?[6] [19, p. 73.]

At this point, M., beginning to lose patience, proceeded to contain Luzin's diagonal saw in a highly eccentric ellipse with tiny minor half-axis $\varepsilon$, and pointed out that as $\varepsilon$ becomes small, the ellipse shrinks down to the diagonal. No room for teeth! Not ready to give up, Luzin responded:

> This is indeed so if $\varepsilon$ is finite, but if $\varepsilon$ is infinitely small... [19, p. 74.]

M.'s "storm of indignation" fell far short of what would one day become the post-revolutionary phlogiston/biped rhetoric:

> I am talking to you for half an hour about limits and not about your actually infinitely small which don't exist in reality. I prove this in my course. Attend it— although for the time being I don't advise you to do so—and you will be convinced of this... [19, p. 74.]

Still unconvinced, Luzin launched into a long soliloquy about filling a cone with gypsum (plaster used for casts), about mathematical idealisation of chemical processes, and how, after removing the cast from the cone,

> we find out that we have not a cone but rather a solid, actually infinitely small... [19, p. 73.]

This was to be Luzin's final comment in that particular conversation. M.'s last suggestion, before stalking away, was that Luzin

> should bring him a jar of that kind of gypsum. [19, p. 76.]

The remarkable conclusion of this exchange occurred some years later, when student Luzin attended a meeting of the mathematical society on Pfaff equations and sat at the back, unobserved by the professors, with M. sitting further forward. As the speaker adroitly manipulated numerous quantities of the form $dx$, $\delta f$, etc. on the blackboard, M., unaware of Luzin's presence, remarked to his neighbor:

> I have always thought that the symbols for exact differentials are special symbols. Look at how he works with them! In his hands they are simply constant numbers: he adds, subtracts, multiplies, substitutes, and transforms them. One can completely forget their origin and operate with them as if they were constant infinitely smalls. [19, p. 77.]

He continued, saying:

---

[6] An allusion to Weierstrass's example of a continuous but nowhere diferentiable function

Copy printed at 1:35 PM on October 25, 2011            5

> it is not at all a hopeless attempt, in the spirit of Hilbert,[7] to axiomatically...
> [19, p. 77.]

but stopped when the speaker was disturbed by the conversation. At this point, Luzin recalled in his letter to Vygodskiĭ that, "A storm erupted in my mind":

> So that's what it is! They teach us, kids, one thing, and they, the grown ups, talk differently to one another. This means that, in fact, to judge by their conversations, things are not so absolutely determined. [19, p. 78.]

With hindsight we know they are not. Luzin continues:

> I looked at them with blazing eyes. I don't know what happened, maybe my stool squeaked … M. suddenly turned around, saw my blazing stare, leaned towards [his neighbor] and said something to him in a low voice. The latter replied in an equally low voice and they [both] fell silent. [19, p. 78.]

As a professional mathematician, Luzin fully understood the need for formalism, but contrasted this with the need for understanding:

> I look at the burning question of the foundations of infinitesimal analysis without sorrow, anger, or irritation. What Weierstrass-Cantor did was very good. That's the way it had to be done. But whether this corresponds to what is in the depths of our consciousness is a very different question. I cannot but see a stark contradiction between the intuitively clear fundamental formulas of the integral calculus and the incomparably artificial and complex work of their "justification" and their "proofs". [19, p. 80.]

We will return to the dialogue between Luzin and Professor M. in Section 10.

## 6. Human thought processes and infinitesimals

Luzin clearly identifies the schism between infinitesimals that seem to make intuitive sense, on the one hand, and the formal definition of limit that gives a sound basis for mathematical analysis, on the other. Once the real numbers have been formally constructed as a complete ordered field, it can be proved that there is no room for infinitesimals in the real number system, so their use was widely condemned. And yet ideas of arbitrarily small quantities continue to be useful in thinking about the calculus because they arise from the natural way in which the brain thinks about variables that become arbitrarily small.

   Mathematical thinking takes place in the human brain where signals take a few milliseconds to pass between neurons to build up a mental conception. Depending on the connections made it takes around a fortieth of a second to see an object and to recognize it. This process continues in time, and we are able to connect together our perceptions and actions as they change dynamically. It happens naturally when drawing a graph with a continuous stroke of a pencil, or looking along a graph to see its changing slope ([32],

---

[7] Hilbert exploited non-Archimedean extension of the reals at about the turn of the century in proving the independence of his axioms of geometry; M. was apparently aware of such a development.

[7]).[8] The natural concept of continuity emerges as a dynamic sense of movement over an interval of time and space and certainly not as a formal definition of a limit at a point. In the same way, when we consider a potentially infinite sequence of values

$$x_1, x_2, \ldots, x_n, \ldots$$

it is natural for us to imagine not just the distinct numerical values, but to think of the n-th term $x_n$ as a dynamically changing entity. Empirical evidence shows how both learners and expert mathematicians imagine such a variable entity to be 'arbitrarily small' (see Cornu [5]). Infinitesimal concepts therefore arise naturally in human thought, causing a conflict between the natural thought processes of learners and the formal modes of proof of mathematical analysis.

While mathematicians may learn to share their formal approach and use it with great success, the transition from intuitive mathematics full of imaginative ideas to formal mathematics based on formal definitions and logical step-by-step deduction presents significant difficulties for many learners (see Pinto & Tall [23], Weber [35]). Fruitful interplay between infinitesimals and student intuitions in the context of 0.999… is explored by R. Ely [8] and Katz & Katz [11], [12]. Related issues in the history and philosophy of mathematics are explored by P. Błaszczyk [1], A. Borovik and M. Katz [2], P. Giordano [9], [10], T. Mormann [22], D. Sherry [28], and Katz & Katz [13], [14], [15].

## 7. A new synthesis of intuition and formalism

Foundational disputes among mathematicians are frequently formulated in purely mathematical terms. However, mathematicians involved in such disputes are not always effective in addressing the transition from intuition to rigour that may be so difficult for learners. Yet the formal approach to mathematics formulated by Hilbert does not ask what the structures are, only what their properties are and what can be deduced from these properties. From this viewpoint, what matters is not what infinitesimals are, but how they behave. An infinitesimal $\varepsilon$ is a (non-zero) element of an ordered field $K$ where $-r < \varepsilon < r$ for all positive rational numbers $r$. If $K$ is an ordered extension field of the real field $\mathbb{R}$, then an infinitesimal will satisfy $-r < \varepsilon < r$ for all positive real numbers $r$.

Infinitesimals cannot fit into the real number system itself, for if $\varepsilon$ is a quantity where $0 < \varepsilon < r$ for all positive real numbers $r$, then $\varepsilon$ cannot be real, for then $r = \tfrac{1}{2}\varepsilon$ is also real and smaller than $\varepsilon$. However, this does not rule out the possibility that an infinitesimal may be an element of an ordered field $K$ which is an ordered extension of $\mathbb{R}$.

---

[8] This analysis by Tall and Katz, developed in [32], is built on Donald's notion of three levels of consciousness in [7].

In this case, any element *k* in an ordered extension field *K* of $\mathbb{R}$ is either infinite (meaning $k > r$ for all $r \in \mathbb{R}$, or else $k < r$ for all $r \in \mathbb{R}$), or it is finite (meaning it lies between two real numbers $a < k < b$). It is straightforward to prove that a finite element *k* is precisely of the form $c + \varepsilon$ where c is real and $\varepsilon$ is either zero or infinitesimal.[9]

For any finite element *k*, its 'standard part' is by definition the unique $c \in \mathbb{R}$, written $c = \mathrm{st}(k)$, such that $k = c + \varepsilon$ where $\varepsilon$ is infinitesimal or zero.

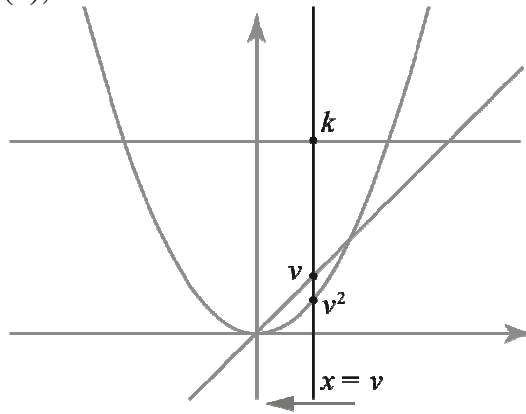

*Figure 2. A helpful sliding vertical line*

As an example, consider the field $\mathbb{R}(x)$ of rational functions in an indeterminate *x*. This field is an extension of the real numbers. The field can be ordered, by defining a rational function *f* to be 'positive' if it is positive in some open interval (0, *a*) for some positive real number *a*. A rational function, in this sense, is either zero, 'positive', or else $-f$ is 'positive'.

The ordered field has a visual representation as graphs in the plane. For any positive real number *k*, the rational functions $y = k$, $y = x$, $y = x^2$ are ordered in the relation $0 < x^2 < x < k$.

By drawing the vertical line $x = v$, the three rational functions meet the line in three points $k$, $v$, $v^2$, where the point $k$ is constant as $v$ varies, but $v$ and $v^2$ are variable, see Figure 2.

A further representation can be obtained by imagining the field $\mathbb{R}(\varepsilon)$ as points on a number line. Clearly infinitesimals are not visible to the naked eye. However, the map $m : \mathbb{R}(\varepsilon) \to \mathbb{R}(\varepsilon)$ given by

$$m(x) = \frac{x - c}{\varepsilon}$$

maps *c* to 0 and $c + \varepsilon$ to 1, thus separating out the images of *c* and $c + \varepsilon$. Following this map by taking the standard part of the image (whenever the image is finite), we obtain $\mathrm{st}(m(x)) \in \mathbb{R}$ (Figure 3).

---

[9] We provide a brief proof. Let $S = \{x \in \mathbb{R} \mid x < k\}$, then *S* is nonempty since it contains *a*, and is bounded above by $b \in \mathbb{R}$, so *S* has a least upper bound *c*, and then the difference $\varepsilon = k - c$ can be shown to be infinitesimal.

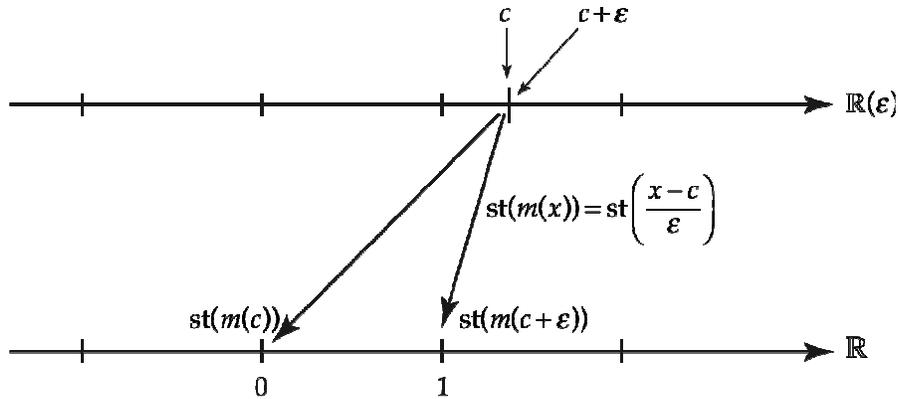

*Figure 3. Visualizing points that differ by an infinitesimal: resolving infinite closeness*

This gives a map to the real numbers which distinguishes the images of the points c and $c + \varepsilon$.

If we imagine such a line as an enhancement of the vertical *y*-axis, then this allows the axis to be imagined as a vertical line with all the elements of $\mathbb{R}(x)$ placed upon it, where now $(0, c + \varepsilon)$ is a fixed point at an infinitesimal distance $\varepsilon$ from the real number $(0, c)$. The point $(0, c + \varepsilon)$ is indistinguishable from $(0, c)$ to the human eye, but it can be distinguished by magnifying the line using the map $m : \mathbb{R}(\varepsilon) \to \mathbb{R}(\varepsilon)$ on the second coordinate and taking the standard part to see a real picture. In this way, we can now imagine the vertical *y*-axis to be a line with fixed infinitesimals that, in a thought experiment, represent the 'final' position of the variable points on the vertical line $x = v$.

This gives four isomorphic representations of an ordered field consisting of rational functions in a single element:

(a) The symbolic system $\mathbb{R}(x)$ of rational functions in an indeterminate *x*;

(b) The graphical system $\mathbb{R}(x)$ of graphs of rational functions in a variable *x*;

(c) The system $\mathbb{R}(v)$ of points on a line where some are constants (the real numbers) and other quantities are variable points determined by where a rational function in *x* meets the line $x = v$. These include infinitesimals such as $v$ and $v^2$, and infinite elements such as $1/v$.

(d) The elements in the ordered field $\mathbb{R}(\varepsilon)$ where $\varepsilon$ is an infinitesimal, which may be represented as fixed points on an extended number line that may be revealed by an appropriate magnification.

In formal terms, all these systems have isomorphic structures that represent the same underlying axiomatic structure: rational functions in a single indeterminate (or variable) with real coefficients. As representations they have very different meanings, for instance, (c) has 'variable' infinitesimals and (d) has fixed infinitesimals, but as formal structures they are isomorphic. To favour one over another is a matter of choice rather than a matter of the underlying formal structure.

## 8. Extending functions beyond the real numbers

To be able to perform the process of differentiation using infinitesimal increments, one needs to be able to extend a given real function *f* from its definition at a real value *x* to a nearby value $x + \varepsilon$ where $\varepsilon$ is an infinitesimal. One can then attempt to define the derivative as

$$f'(x) = \text{st}\left(\frac{f(x+\varepsilon) - f(x)}{\varepsilon}\right)$$

For example, if $f(x) = x^2$, then $f'(x) = \text{st}(2x + \varepsilon) = 2x$.

If *f* is a rational function, it is easy to substitute $x + \varepsilon$ for *x* in the extension field $\mathbb{R}(\varepsilon)$ and compute the value of the function using algebra. However, more general functions such as $\sin x$ or $e^x$ cannot be extended in this field. Meanwhile, both extensions are possible in the field of series $\mathbb{R}((\varepsilon))$ in an infinitesimal $\varepsilon$, which allow a finite number of terms in $1/\varepsilon$, of the form

$$a_{-N}\varepsilon^{-N} + \ldots + a_{-1}\varepsilon^{-1} + a_0 + a_1\varepsilon + \ldots + a_n\varepsilon^n + \ldots$$

These form an ordered field in which $\varepsilon$ is an infinitesimal and can be used to extend any function expressible as a power series to a function $f(x + \varepsilon)$. This field was called the *superreals* by Tall [30]. This is sufficient to deal with analytic functions (given by power series), which is essentially strong enough for combinations of standard functions in the calculus, but not for general functions in mathematical analysis.

The more general problem is to extend every real function *f* to take on values in an extension field with infinitesimals.

We could begin by working with all sequences of real numbers $(x_n)$ and then operate with them term by term. Then a sequence such as (1, 2, 3, ...) might be considered as an infinite number and its inverse would be $(1, \frac{1}{2}, \frac{1}{3}, \ldots)$ which would be an infinitesimal. The constant sequence (*k*, *k*, *k*, ...) for any $k \in \mathbb{R}$ could then be identified with *k*, to let us embed $\mathbb{R}$ in the set of all such sequences. The problem is that such a system does not operate as an ordered field. For instance, even though we we might like to think of the sequence (1, 2, 3, ...) as being bigger than any real number *k*, its initial terms might be less than *k* and the *n*th term would only exceed *k* once $n > k$.

The first step towards equivalence would be to say that two sequences $(a_n)$, $(b_n)$ are equivalent if they are equal for all but a finite number of terms. If we denoted the equivalence class containing $(a_n)$ by $[a_n]$, then we would have $[a_n] = [b_n]$ if and only if $a_n = b_n$ for all but a finite number of *n*.

This would give us a surprisingly good beginning to the problem, for if we let $\omega$ be the equivalence class of (1, 2, 3, ...), then we would have $\omega > k$ for any real number $k$ because all but a finite number of its terms are greater than $k$ (by identifying $k$ with the sequence $(k, k, k, …)$).

Furthermore $\omega + 1$ would be (2, 3, ...), giving a situation in which every term of (2, 3, ...) is 1 bigger than the corresponding term of (1, 2, ...), so $\omega + 1 > \omega$, unlike cardinal infinities, where $\aleph_0 + 1 = \aleph_0$. The term $\omega^2$ would be far bigger, $1/\omega$ (the sequence $(1, \frac{1}{2}, \frac{1}{3}, …)$) would be infinitesimal and $1/\omega^2$ would be a smaller infinitesimal still. We would even have a natural extension of a set $D$ to the larger set *$D$ consisting of all equivalence classes $[x_n]$ where $x_n \in \mathbb{R}$ and any function $f : D \to \mathbb{R}$ can be extended to $f : *D \to *\mathbb{R}$ by defining f($[x_n]$) to be the equivalence class $[f(x_n)]$.

We still need to do more. As well as dealing with sequences that nicely tend to a limit, we need to deal with *every* sequence of real numbers, no matter how it is defined. For instance, the sequence (0, 1, 0, 1, 0,...) equals 0 on the set $O$ of odd numbers, but equals 1 on the set $E$ of even numbers. We need to assign it to an appropriate equivalence class. If we make the decision focusing on the odd numbers, it will be equivalent to 0, but if we make a decision focusing on the even numbers, it will be equivalent to 1. The consequence is that to define the equivalence relation fully, we must make a *choice*.

Making such a choice may seem strange at first, but it is only a technical device to make a decision in all cases so that we decide that every $[x_n]$ represents a specific finite or infinite quantity. Since the set of terms in the sequence $(x_n)$ is infinite, then at least one of the following three possibilities must hold:

(i) there is an infinite subsequence tending to $-\infty$.

(ii) there is an infinite subsequence tending to $+\infty$.

(iii) for some $A, B \in \mathbb{R}$, where $A < B$, there is an infinite number of terms between $A$ and $B$.

All three possibilities may occur, as in the case of the sequence $(a_n)$ given by $(-1, 2, \frac{1}{3}, -4, 5, \frac{1}{6}, …)$, where $a_n = -n$ for $n = 3N - 2$, $a_n = n$ for $n = 3N - 1$, and $a_n = \frac{1}{n}$ for $n = 3N$, as $N$ increases through 1, 2, ... .

In general, if case (i) holds, then we may choose $[a_n]$ to be negative infinite by taking the decision set to be the set of $n$ for which a subsequence of terms tends to $-\infty$. (In the example, the decision set {1, 4, 7, ...} gives $[a_n] = -\omega$. In case (ii), we may choose $[a_n]$ to be positive infinite. (In the example, the decision set {2, 5, 8, ...} gives $[a_n] = \omega$. In case (iii), the terms have a subsequence tending to a finite number $L$ where $A \leq L \leq B$ and we may choose to make the decision on the set related to this subsequence, which gives $[a_n]$ equal to $L$ plus an infinitesimal. (In the example, the decision set is {3, 6, 9, ...}, $L = 0$ and $[a_n]$ is the infinitesimal $1/\omega$.)

The major problem is to choose all the decision sets in such a way that all the choices can be made in a consistent manner.

## 9. Making a serious choice

To make a coherent decision in all possible cases requires us to formulate a full collection of decision sets and say that a property $P(x)$ is true for a sequence $x = (x_n)$ if $P(x_n)$ is true for all $n$ in a particular decision set. If $S$ is chosen to be a decision set, we will say that $S$ is *decisive*.

First, we stipulate that a finite set cannot be decisive, while every cofinite set $S$ (a set consisting of the whole of $\mathbb{N}$ except for a finite number of elements) is necessarily decisive:

(0) If $S$ is finite, then $S$ is not decisive.

(1) If $\mathbb{N} \setminus S$ is finite, then $S$ is decisive.

Next, if we decide that $P(x)$ is true because $P(x_n)$ is true for all $n$ in a decision set $S$, then it may also be true in a larger set $T$. For the sake of coherence, $T$ should also be chosen to be decisive:

(2) If $S$ is decisive and $S \subset T \subset \mathbb{N}$, then $T$ is also decisive.

We further require that

(3) If $S$ and $T$ are decisive then $S \cap T$ is decisive.

Thus, the intersection of a pair of decisive sets should be decisive, as well. A collection of sets satisfying (0)–(3) is a *filter* on $\mathbb{N}$. It is relatively easy to construct a filter step by step. Just start with property (1) to include all subsets of $\mathbb{N}$ whose complements are finite. Then, if any new set $U$ is added, so must all sets of natural numbers that contain $U$, and any intersection of $U$ with a subset already in the filter. For example, if we start with the sets required by (1) and add the set $E$ of all even numbers, then we need to include any set containing the even numbers, as well as any subset of these sets formed by omitting a finite number of elements. The new set of sets now satisfies (0)–(3) and so it is again a filter.

The serious problem comes with expanding such a filter to satisfy the following additional requirement:

(4) For each subset $S$ of $\mathbb{N}$, one of the two sets $S$ and $\mathbb{N}\setminus S$ must be decisive.

This requires an infinite number of decisions to be made and seems impossible for a human being to accomplish in a finite lifetime. But that does not mean that we cannot imagine it happening.

A filter satisfying properties (1)–(4) is called a (nonprincipal) *ultrafilter*.

(Axiom (4) now renders (0) redundant as it follows from a combination of (1) and (4).)

The existence of such an ultrafilter is guaranteed by the axiom of choice, see Tarski [33]. The choice is not unique. For example, if we choose the odd numbers to be decisive, then for $a = [1, 0, 1, 0, ...]$ and $b = [0, 1, 0, 1, ...]$ we have $a > b$ and, if not, then by condition (4), the set of even numbers is decisive and we have $b > a$.

The set of equivalence classes is denoted by $^*\mathbb{R}$ and is called a field of hyperreal numbers.[10] The properties (1)–(4) guarantee properties that might be expected. For instance, comparing the element $\omega = [1, 2, \ldots, n, \ldots]$ with $k = [k, k, \ldots, k, \ldots]$ gives $n > k$ for all but a finite number of $n$, so by property (1), we have $\omega > k$ for all real numbers $k$ and hence $\omega$ is infinite. Similarly $1/\omega$ is a positive infinitesimal because $0 < 1/n < k$ for all $n > 1/k$ for any positive real number $k$.

Such a system of hyperreals together with the extension of any subset $D$ to the subset

$$^*D = \{[x_n] \in {^*\mathbb{R}} \mid x_n \in D\}$$

gives an extended map $f : {^*D} \to {^*\mathbb{R}}$ defined by $f([x_n]) = [f(x_n)]$.

## 10. Looking closely at Luzin's saw

Returning to the earlier difference of opinion between Luzin and Professor M., we now see that M. is correct in declaring that infinitesimals cannot be perceived nor made out of gypsum at their original size. However, the subsequent development of mathematical theory shows that he rejects the possibility of allowing students to give infinitesimals a mathematical meaning that can be perceived not only in imagination but also in a physical picture.

Let us imagine a finite sawtooth with $n$ equal steps from $(0, 0)$ to $(1, 1)$ as in Figure 4.

---

[10] We say *a* field of hyperreals, rather than *the* field of hyperreals, because different choices are possible. If one assumes that the continuum hypothesis is true, then our ultrapower construction produces a unique ordered field up to isomorphism.

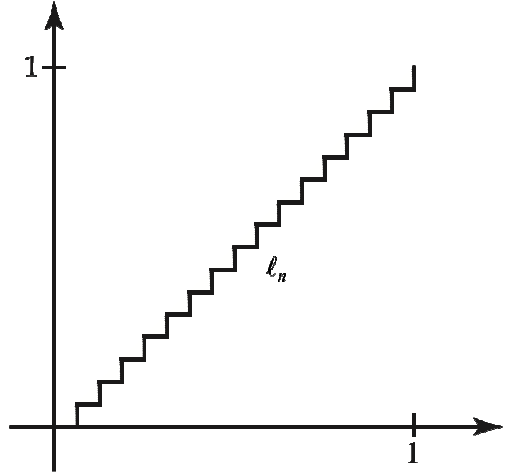

*Figure 4. Luzin's finite saw, parametrised as $\ell_n : [0,1] \to \mathbb{R}^2$*

This can be parametrised as

$$\ell_n : [0,1] \to \mathbb{R}^2$$

by tracing along it with a finger in time $t$ from 0 to 1 in which, for $k = 0, ..., n-1$, the $k$th step begins at the point $(t, t)$ at $t = k/n$, moves up to $(t, t + \frac{1}{n})$ and then to $(t + \frac{1}{n}, t + \frac{1}{n})$ at time $t + \frac{1}{n}$.

The graph of $\ell_n(t)$ for $0 \leq t \leq 1$ is the $n$th finite Luzin saw drawn in the coordinate plane. His saw with infinitesimal teeth may be conceived as $\ell_N$ for infinite $N \in {}^*\mathbb{N}$. This is quite natural in ${}^*\mathbb{R}^2$, but as Luzin envisaged it, he thought of $N$ as being a *countable* set, and $N$ cannot be an infinite cardinal, because an infinite cardinal does not have an inverse $1/N$. However, this does not mean he cannot *imagine* infinitesimal quantities in an extension field, only that there is still work to be done to formulate a formally coherent structure (such as the hyperreals).

Professor M. is correct in asserting that the limit of $(\ell_n)$ is the straight line joining (0, 0) to (1, 1). However, the graph of $\ell_N$ in the extended ${}^*\mathbb{R}^2$-plane has an infinite number $N$ of saw-teeth. The $k$th sawtooth (where $k$ may now be finite or infinite) starts at $(t, t)$ for $t = k/N$, moves up to $(t, t + \frac{1}{N})$ and then to $(t + \frac{1}{N}, t + \frac{1}{N})$ at time $t + \frac{1}{N}$.

When we magnify by the factor $N$ pointing at $(X, Y)$ where $X = k/N$, $Y = k/N$, we get a map $m : {}^*\mathbb{R}^2 \to {}^*\mathbb{R}^2$ in the form

$$m(x, y) = (N(x - X), N(y - Y)).$$

which gives

$$m(X, Y) = (0,0), \text{ and } m(X + 1/N, Y + 1/N) = (1,1)$$

so that the sides of the saw-tooth of size $1/N$ are magnified to unit lengths, giving the picture in Figure 5.

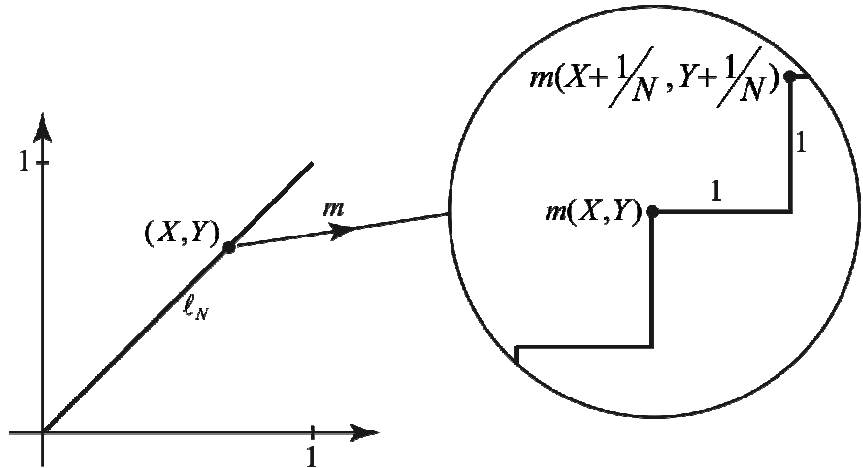

*Figure 5. Luzin's infinitesimal saw*

This reveals that the original line seen through perceptual human eyes (without magnification) does indeed look like M.'s diagonal line and, from M.'s viewpoint seeing only points in $\mathbb{R}$, the limit of the finite saw-teeth $(\ell_n)$ is indeed just the diagonal line. In the 'real' world of Cantor, with only the complete ordered field $\mathbb{R}$ Luzin would have to work a little harder to give an everywhere continuous nowhere differentiable curve. Rather than work with the Luzin sawteeth along the diagonal, let us look look at a sawtooth along the real line to get the first saw-tooth

$$s_1(x) = \begin{cases} x & \text{for } 0 \leq x \leq \tfrac{1}{2} \\ 1-x & \text{for } \tfrac{1}{2} \leq x \leq 1 \end{cases}, \quad s(n+x) = s(x) \text{ for } n \in \mathbb{Z}.$$

And then define successive saw-teeth as

$$s_n(x) = s_1(2^{n-1} x) / 2^{n-1}.$$

These give successive teeth half the size of each previous one (figure 6). This variant of Luzin's finite saw-teeth has limit zero as *n* increases:

$$\lim_{n \to \infty} s_n(x) = 0,$$

to give a straight line, just as M. declared it would be (Figure 6).

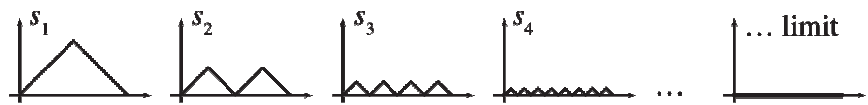

*Figure 6: The real limit of Luzin's sawteeth using M.'s mathematical analysis*

However, if we *add together* the saw-teeth to get

$$\text{bl}(x) = \lim_{n \to \infty} \sum_{k=1}^{n} s_n(x),$$

then we get the blancmange function [31], identified by Takagi [29] in 1903, re-invented and generalized by van der Waerden [34] in 1930. This function is everywhere continuous and nowhere differentiable (Figure 7).

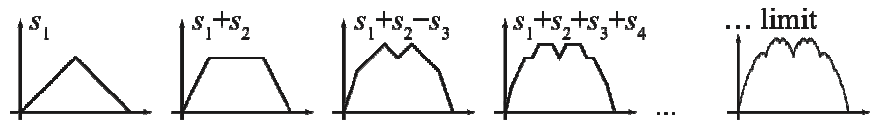

*Figure 7: The real limit adding Luzin's sawteeth
to give a real continuous, non-differentiable function*

The moral here is that the context in which one is working affects the nature of the mathematics. If one works in a context that only allows a complete ordered field $\mathbb{R}$, then M. is totally vindicated, along with Cantor and the mathematical culture of standard analysis. However, in the broader context of formal mathematics in extension fields of $\mathbb{R}$, infinitesimals *must* occur. Here the intuition of Luzin can be formally defined in an extended system (that he did not have at the time) in which the graph of $|_N$ has infinitesimal steps as he suggested. Furthermore, these steps can be represented in a physical drawing using a formally defined magnification, for all to see.

## 11. Conclusion

The approach outlined here allows the handling of sequences tending to zero or to infinity to be conceived in terms of infinitesimal and infinite quantities in a formal manner that is consonant with intuitive ideas of infinitesimals and infinite numbers in the calculus. The construction requires standard set theory together with the axiom of choice. Adjoining the axiom of choice does not introduce any contradictions in the sense that if standard set theory is consistent, so is the system when the axiom of choice is added.

    While Professor M. and his modern counterparts believe that they are protecting the purity of mathematics by telling students that 'the actual infinity does not exist', this denies their students' right to imagine infinitesimal and infinite quantities in their mind's eye and to link their intuitive vision at some stage to a full formal approach.

    Ideas of infinitesimals being generated by 'variables that tend to zero' were introduced by Cauchy and offer a meaning that is still used in applications today. Likewise 'variables that can be arbitrarily small' evoke a natural sense of dynamic limit processes appropriate for the calculus that can later be transformed either to standard arguments in mathematical analysis or infinitesimal methods based on the hyperreals. Meanwhile, as we saw in the structure theorem for any ordered extension field *K* of $\mathbb{R}$, any such field can be imagined as an enriched number line with fixed infinitesimals that can be distinguished using a visual picture of a formal magnification. Moreover, the conceptions of infinitesimals as variables or as fixed quantities are different representations of the same underlying concept because *any* ordered field

extension can be visualized as an enriched number line and magnified to see its infinitesimal quantities.

Most modern mathematicians now admit the axiom of choice, in the knowledge that it offers theoretical power without introducing contradictions that did not exist before. Is it not time to allow infinitesimal conceptions to be acknowledged in their rightful place, both in our fertile mathematical imagination and in the power of formal mathematics, enriched by the axiom of choice?

**References**

[1] Błaszczyk, P.: Nonstandard analysis from a philosophical point of view. In Non-Classical Mathematics 2009 (Hejnice, 18-22 june 2009), pp. 21-24.

[2] Borovik, A.; Katz, M.: Who gave you the Cauchy--Weierstrass tale? The dual history of rigorous calculus. Foundations of Science, 2011, see http://dx.doi.org/10.1007/s10699-011-9235-x and http://arxiv.org/abs/1108.2885

[4] Cauchy, A. L.: Cours d'Analyse de l'Ecole Royale Polytechnique. Paris: Imprimérie Royale, 1821 (reissued by Cambridge University Press, 2009.)

[5] Cornu, B.: Limits, pp. 153-166, in *Advanced mathematical thinking*. Edited by D. O. Tall. Mathematics Education Library, 11. Kluwer Academic Publishers Group, Dordrecht, 1991.

[6] Courant, R.: *Differential and integral calculus*. Vol. I. Translated from the German by E. J. McShane. Reprint of the second edition (1937). Wiley Classics Library. A Wiley-Interscience Publication. John Wiley & Sons, Inc., New York, 1988.

[7] Donald, M.: *A Mind So Rare*. New York: Norton & Co., 2001.

[8] Ely, R.: Nonstandard student conceptions about infinitesimals. Journal for Research in Mathematics Education 41 (2010), no. 2, 117-146.

[9] Giordano, P.: Infinitesimals without logic. Russ. J. Math. Phys. 17 (2010), no. 2, 159-191.

[10] Giordano, P.; Katz, M.: Two ways of obtaining infinitesimals by refining Cantor's completion of the reals, 2011, see http://arxiv.org/abs/1109.3553

[11] Katz, K.; Katz, M.: Zooming in on infinitesimal 1-.9.. in a post-triumvirate era. Educational Studies in Mathematics 74 (2010), no. 3, 259-273. See arXiv:1003.1501.

[12] Katz, K.; Katz, M.: When is .999… less than 1? The Montana Mathematics Enthusiast 7 (2010), No. 1, 3--30.

[13] Katz, K.; Katz, M.: A Burgessian critique of nominalistic tendencies in contemporary mathematics and its historiography. Foundations of Science , 2011, see http://dx.doi.org/10.1007/s10699-011-9223-1 and http://arxiv.org/abs/1104.0375

[14] Katz, K; Katz, M.: Cauchy's continuum. Perspectives on Science 19 (2011), no. 4, 426-452. See http://dx.doi.org/10.1162/POSC_a_00047 and http://arxiv.org/abs/1108.4201

[15] Katz, K.; Katz, M.: Stevin numbers and reality. Foundations of Science 2011, see http://dx.doi.org/10.1007/s10699-011-9228-9 and http://arxiv.org/abs/1107.3688

[16] Kiepert, L.: *Grundriss der Differential-und Integralrechnung. I. Teil: Differentialrechnung.* Helwingsche Verlagsbuchhandlung. Hannover, 12th edition. 1912. XX, 863 S.

[17] Klein, F.: *Elementary Mathematics from an Advanced Standpoint. Vol. I. Arithmetic, Algebra, Analysis*. Translation by E. R. Hedrick and C. A. Noble [Macmillan, New York, 1932] from the third German edition [Springer, Berlin, 1924]. Originally published as *Elementarmathematik vom höheren Standpunkte aus*, Leipzig, 1908.

[18] Lakoff, G., Núñez, R.: *Where mathematics comes from. How the embodied mind brings mathematics into being*. Basic Books, New York, 2000.

[19] Luzin, N. N. (1931) Two letters by N. N. Luzin to M. Ya. Vygodskii. With an introduction by S. S. Demidov. Translated from the 1997 Russian original by Shenitzer. *Amer. Math. Monthly* 107 (2000), no. 1, 64–82.

[20] Medvedev, F. A.: *Nonstandard analysis and the history of classical analysis. Patterns in the development of modern mathematics* (Russian), 75–84, "Nauka", Moscow, 1987.

[21] Medvedev, F. A.: Nonstandard analysis and the history of classical analysis. Translated by Abe Shenitzer. *Amer. Math. Monthly* 105 (1998), no. 7, 659 – 664.

[22] Mormann, T.: Idealization in Cassirer's philosophy of mathematics. Philos. Math. (3) 16 (2008), no. 2, 151-181.

[23] Pinto, M. M. F., Tall, D. O.: Student constructions of formal theory: giving and extracting meaning. In O. Zaslavsky (Ed.), *Proceedings of the 23rd Conference of PME*, Haifa, Israel, 4, (1999), 65–73.

[24] Robinson, A.: *Non-standard analysis*. North-Holland Publishing Co., Amsterdam 1966.

[25] Robinson, A., Roquette, P.: On the finiteness theorem of Siegel and Mahler concerning diophantine equations. *J. Number Theory* 7 (1975), 121–176.

[26] Roquette, P.: Numbers and models, standard and nonstandard. *Math Semesterber* 57 (2010), 185–199.

[27] Roquette, P.: Nonstandard aspects of Hilbert's irreducibility theorem. In: Saracino, D.H., Weispfenning, B. (eds.) *Model Theor. Algebra*, Mem. Tribute Abraham Robinson, vol. 498 of Lect. Notes Math. pp. 231–275. Springer, Heidelberg, 1975.

[28] Sherry, D.: The wake of Berkeley's Analyst: *rigor mathematicae*? Stud. Hist. Philos. Sci. 18 (1987), no. 4, 455--480.

[29] Takagi, T.: A simple example of the continuous function without derivative, *Proc. Phys.-Math. Soc,* Tokyo Ser. II 1 (1903), 176–177.

[30] Tall, D. O.: Looking at graphs through infinitesimal microscopes, windows and telescopes, *Mathematical Gazette*, 64 (1980), 22–49.

[31] Tall, D. O.: The blancmange function, continuous everywhere but differentiable nowhere, *Mathematical Gazette*, 66 (1982), 11–22.

[32] Tall, D. O., Katz, M.: A Cognitive Analysis of Cauchy's Conceptions of Continuity, Limit, and Infinitesimal, with implications for teaching the calculus (submitted for publication, 2011). Available on the internet at: http://www.warwick.ac.uk/staff/David.Tall/downloads.html.

[33] Tarski, A., Une contribution à la théorie de la mesure, *Fund. Math.* 15 (1930), 42–50.

[34] Waerden, B. L. van der: Fan einfaches Beispiel einer nichtdifferenzierbaren stetigen Funktion, *Math. Z.* 32 (1930), 474–475.

[35] Weber, K.: Traditional instruction in advanced mathematics courses: a case study of one professor's lectures and proofs in an introductory real analysis course, *Journal of Mathematical Behavior*, 23 (2004), 115–133.